\documentclass[11pt]{amsart}
\usepackage{amssymb}
\usepackage{calrsfs}
\usepackage[all]{xy}
\usepackage{tikz}
\usetikzlibrary{arrows,arrows.meta,positioning}
\usepackage{xcolor}
\usepackage[letterpaper,margin=1in]{geometry}

\usepackage{hyperref}

\newtheorem{dummy}{dummy}[section]

\newtheorem{theorem}[dummy]{Theorem}

\newtheorem{conjecture}[dummy]{Conjecture}

\theoremstyle{definition}

\newtheorem{example}[dummy]{Example}
\newtheorem{remark}[dummy]{Remark}

\newtheorem{expectation}[dummy]{Expectation}

\newcommand{\bC}{\mathbb{C}}

\newcommand{\bL}{\mathbb{L}}

\newcommand{\bP}{\mathbb{P}}

\newcommand{\bT}{\mathbb{T}}

\newcommand{\cA}{\mathcal{A}}

\newcommand{\cC}{\mathcal{C}}

\newcommand{\cE}{\mathcal{E}}

\newcommand{\cG}{\mathcal{G}}

\newcommand{\cL}{\mathcal{L}}

\newcommand{\cO}{\mathcal{O}}

\newcommand{\cS}{\mathcal{S}}

\DeclareMathOperator{\Fun}{Fun}
\DeclareMathOperator{\Sing}{Sing}

\DeclareMathOperator{\Perf}{Perf}
\DeclareMathOperator{\Qcoh}{QCoh}
\DeclareMathOperator{\Spec}{Spec}

\DeclareMathOperator{\Coh}{Coh}

\DeclareMathOperator{\Hom}{Hom}

\DeclareMathOperator{\coker}{coker}

\DeclareMathOperator{\HH}{HH}

\DeclareMathOperator{\Mod}{Mod}

\DeclareMathOperator{\Sym}{Sym}
\DeclareMathOperator{\Indcoh}{IndCoh}

\newcommand{\Fuk}{\mathrm{Fuk}}

\newcommand{\pt}{\mathrm{pt}}

\newcommand{\dgcat}{\mathrm{DGCat}}

\newcommand{\locsys}{\mathrm{LocSys}}
\newcommand{\locsyscat}{\mathrm{LocSysCat}}
\newcommand{\shvcat}{\mathrm{ShvCat}}
\newcommand{\indcohshvcat}{\mathrm{IndCohShvCat}}
\newcommand{\im}{\mathrm{im}}

\newcommand{\nFuk}[1]{{#1}\mathrm{Fuk}}
\newcommand{\Symp}{\mathrm{Symp}}
\newcommand{\Ham}{\mathrm{Ham}}
\newcommand{\nSymp}[1]{{#1}\mathrm{Symp}}
\newcommand{\frakg}{\mathfrak{g}}
\newcommand{\crit}{\mathrm{Crit}}
\newcommand{\critv}{\mathrm{Critv}}
\newcommand{\schobers}{\mathrm{Schobers}}
\newcommand{\nHom}[1]{{#1}\mathrm{Hom}}
\newcommand{\bone}{\mathbf{1}}

\newcommand{\rB}{\mathrm{B}}
\newcommand{\rC}{\mathrm{C}}
\newcommand{\rH}{\mathrm{H}}
\newcommand{\rT}{\mathrm{T}}
\newcommand{\QH}{\mathrm{QH}}
\newcommand{\SH}{\mathrm{SH}}
\newcommand{\HF}{\mathrm{HF}}
\newcommand{\BFM}{\mathrm{BFM}}
\newcommand{\KRS}{\mathrm{KRS}}
\newcommand{\Prl}{\mathrm{Pr}^{\mathrm{L}}}

\usepackage{setspace}
\title[Higher Fukaya categories]{Speculations on higher Fukaya categories}

\author{James Pascaleff}
\address{James Pascaleff, Department of 
Mathematics, University of Illinois at Urbana-Champaign, IL, US}
\email{jpascale@illinois.edu}

\author{Nicol\`o Sibilla}
\address{Nicol\`o Sibilla, 
SISSA,
Via Bonomea 265, 
34136 Trieste, Italy}
\email{nsibilla@sissa.it}

\begin{document}

\begin{abstract}
  We investigate a possible theory of higher Fukaya categories associated to $n$-shifted symplectic stacks, where $n \geq 0$. We consider two paradigmatic cases, the shifted cotangent stack of a smooth manifold and the coadjoint stack of a compact Lie group, drawing connections to the work of Teleman and 3D mirror symmetry. Our evidence includes some new results in $1$-shifted symplectic geometry.
\end{abstract}

\maketitle

\section{Introduction}

The Fukaya category $\Fuk(X,\omega)$ of a symplectic manifold is widely accepted as the mathematical formalization of the ``category of boundary conditions in the 2-dimensional A-model'' with target $(X,\omega)$. As such, this category is one of the sides of the homological mirror symmetry correspondence (2D HMS). In recent years, there has been much progress in understanding three-dimensional mirror symmetry, including its homological version (3D HMS), where the ``2-category of boundary conditions in the 3-dimensional A-model'' plays a central role \cite{bullimore-et-al,webster-yoo,gammage-hilburn-mazel-gee}. To our knowledge, there does not seem to be a definition of this 2-category that has achieved general acceptance. 

Motivated by these questions, we have sought natural constructions of higher categories that are ``of Fukaya type,'' meaning that their definition is fundamentally based on the intersection theory of Lagrangians. It turns out that, thanks to the development of shifted symplectic geometry \cite{ptvv}, there are a wealth of constructions that can be used for this purpose. However, making our proposals rigorous will still require more foundational work, for reasons that will become clear. Thus we must confine ourselves to a tentative sketching of the theory, but one that leads to a rather neat picture. 

The higher Fukaya categories that we shall consider have for their target space an $(n-1)$-shifted symplectic stack $X$. We should note that, for our desired applications this stack $X$ should be taken in the $C^{\infty}$ category, rather than the algebraic or holomorphic one (see \cite{steffens}). The output is a $(\infty,n)$-category $\nFuk{n}(X)$:
\begin{equation*}
  \text{$(n-1)$-shifted symplectic stack $X$} \rightsquigarrow \text{$(\infty,n)$-category $\nFuk{n}(X)$.}
\end{equation*}
Note that the category number is one more than the shift. The objects of this category are $(n-1)$-shifted Lagrangians $f : L \to X$ (with extra data), and their higher morphisms are computed by iteratively taking derived intersections, with a crucial exception: once we reach the level of $0$-shifted Lagrangians in $0$-shifted symplectic manifold, we use Floer theory to produce a vector space (rather than a $(-1)$-shifted symplectic stack).

To obtain a $2$-category, $X$ should be a $1$-shifted symplectic stack. The simplest such object is $\pt[1]$, the $1$-shifted point. By \cite{ptvv}, $1$-shifted Lagrangians in $\pt[1]$ are nothing but $0$-shifted symplectic stacks. Thus, we are led to expect
\begin{equation*}
  \nFuk{2}(\pt[1]) \supseteq \Symp,
\end{equation*}
where $\Symp$ is the Weinstein $(\infty,2)$-category whose objects are symplectic manifolds, and where the morphisms and compositions are defined using the quilts of Wehrheim-Woodward \cite[and subsequent papers]{wehrheim-woodward}; see \cite{abouzaid-bottman} for a review. Going back to the general case, we can view $\nFuk{2}(X)$ as a relative version of $\Symp$; this in particular implies that the technical challenges involved in defining $\nFuk{2}(X)$ cannot be less than those involved in defining $\Symp$.

Having some idea of what $\nFuk{2}(X)$ might look like, we can return to the motivating question: does this category have anything to do with the category of boundary conditions for the $3$-dimensional A-model? At first glance, it appears to be different. First of all, in many works on 3D HMS, the target for the 3D A-model is a holomorphic symplectic or hyperk\"{a}hler manifold. Thus the target is expected to be $0$-shifted symplectic in the holomorphic or algebraic category, so the shift is different! Floer-theoretic proposals for this kind of category has been advanced by Bousseau \cite{bousseau}, Doan-Rezchikov \cite{doan-rezchikov}, and Kontsevich-Soibelman \cite{kontsevich-soibelman}.

The work of Ki Fung Chan and Conan Leung \cite{chan2024syzmirrorsnonabelian3d, chan20243dmirrorsymmetrymirror} proposes an approach that has similar features to our own, in the sense that they also seek to understand 3D mirror symmetry recursively in terms of 2D mirror symmetry.

The work of Elliott-Yoo \cite{elliott-yoo-langlands,elliott-yoo-sing} on Geometric Langlands is related to 3D mirror symmetry, and it also makes extensive use of shifted symplectic geometry. The relationship of their perspective to the present one is not entirely clear as of this writing, but it should be pursued in the future.

On the other side, there are natural candidates for the $2$-category of boundary conditions in the $3$-dimensional B-model. Teleman \cite{teleman-icm} uses the $2$-category of Kapustin-Rozansky-Saulina \cite{kapustin-rozanksy-saulina} defined for a holomorphic symplectic target, and Gammage--Hilburn--Mazel-Gee \cite{gammage-hilburn-mazel-gee} use a $2$-category of Ind-coherent sheaves of categories. 

The discrepancy between our $1$-shifted targets and the $0$-shifted targets that are more commonly considered could be taken to suggest that, as natural as our categories may be, they are not the right ones for 3D HMS. To the contrary, by analyzing some examples, we discovered that our theory in some cases reproduces structures that are involved in some of the paradigmatic examples of 3D HMS. In this paper we shall present the following examples.
\begin{enumerate}
\item In their work on 3D HMS, Gammage--Hilburn--Mazel-Gee \cite{gammage-hilburn-mazel-gee} used a certain category of \emph{schobers} to represent the 3D A-model. A schober that has no singular points is essentially a \emph{local system of categories}. Our recent work \cite{pps-loc-sys,pps-koszul}, joint with Pavia, is devoted to a thorough exploration of the theory of local systems of higher categories.

  Symplectic geometry supplies many examples of such structures. Given a Hamiltonian fibration $\pi : X \to B$, the Fukaya categories of the fibers form a local system of categories over $B$ \cite{savelyev-global-fuk}. If we allow $\pi : X \to B$ to be a symplectic Lefschetz fibration, we obtain a schober \cite{seidel-fcplt}. We shall find that, in both cases, $\pi$ actually defines an object of
  \begin{equation*}
    \nFuk{2}(\rT^{*}[1]B),
  \end{equation*}
  where $\rT^{*}[1]B$ denotes the $1$-shifted cotangent bundle of $B$. These are novel results in the field of shifted symplectic geometry, and we provide sketches of proofs.

  From this, we conjecture that $\nFuk{2}(\rT^{*}[1]B)$, properly defined, contains the $2$-category of schobers on $B$.

\item Teleman's groundbreaking work \cite{teleman-icm, teleman-coulomb} has inspired much of the subsequent work in this area. He considers the case of Hamiltonian $G$-manifolds $M$. The moment map $\mu : M \to \frakg^{*}$ is $G$-equivariant, and so induces a morphism $[M/G]\to [\frak{g}^{*}/G]$. There is a natural $1$-shifted symplectic structure on $[\frakg^{*}/G]$ making this morphism Lagrangian. Hence the Hamiltonian $G$-manifold $M$ determines an object in
  \begin{equation*}
    \nFuk{2}([\frakg^{*}/G]).    
  \end{equation*}

  One of Teleman's most remarkable discoveries is that the $G$-equivariant quantum cohomology $\QH_{G}^{*}(M)$ carries an action of $\rH_{*}^{G}(\Omega G)$ (also known as functions on the Bezrukavnikov-Finkelberg-Mirkovic space $\BFM(G^{\vee})$ of the Langlands dual group $G^{\vee}$). We shall show that the existence of this action can be deduced from formal higher-categorical manipulations taking place in the $3$-category
  \begin{equation*}
    \nFuk{3}(\pt[2])
  \end{equation*}
  where $\pt[2]$ is the $2$-shifted point. The key here is that $[\frakg^{*}/G]$ is an object of this $3$-category, and the Hamiltonian $G$-space $M$ determines a $1$-morphism.  
\end{enumerate}

In our view, these examples favor the claim that the higher Fukaya category has a natural place in the story of 3D HMS. In Section \ref{sec:coadjoint}, we shall argue for a direct connection between the $2$-category category $\nFuk{2}([\frakg^{*}/G])$ and the 3D B-model $2$-category of $\BFM(G^{\vee})$, analogous to the \emph{intrinsic mirror symmetry} philosophy of Gross-Siebert \cite{gross-siebert}.

These examples also indicate that the higher Fukaya category can potentially unify the theory of Lefschetz fibrations and schobers with Teleman's perspective on gauge theory and mirror symmetry.

\begin{remark}
  We reiterate that our constructions are intended to take place in the setting of derived higher stacks over $C^{\infty}$-manifolds, not in the more familiar algebraic setting of \cite{ptvv} and most related works. See \cite{steffens} for some initial foundations of the $C^{\infty}$ theory. We shall take for granted that the core features of derived symplectic geometry carry over to the $C^{\infty}$ setting, allowing us to quote results from \cite{ptvv} and so on, but strictly speaking these results are merely analogous to the ones we need.
\end{remark}

\subsection*{Acknowledgements}

The authors would like to thank the organizers of the workshop ``Mathematics inspired by String Theory,'' March 2025, at the Chinese University of Hong Kong and the Hong Kong Geometry Colloquium where this work was presented. They specifically thank Conan Leung, Michael McBreen, and Philsang Yoo for their interest in this work.

The authors also thank Emanuele Pavia for his collaboration on the topic of higher-categorical local systems, which are a facet of the picture presented here.

JP was partially support by the Simons Foundation.

\section{What is a ``Fukaya-type'' category?}

Our answer to the question of the title is that a Fukaya-type category should have Lagrangians as objects, and the morphisms should be defined using Lagrangian intersection theory. Of course, there are multiple ways that the phrase ``Lagrangian intersection theory'' could be interpreted. Among them are:
\begin{enumerate}
\item Lagrangian Floer cohomology. Here $L_{0}$ and $L_{1}$ are Lagrangians in an ordinary symplectic manifold $X$ (with extra data), and the Floer cohomology $\HF^{*}(L_{0},L_{1})$ is ``semi-infinite Morse homology'' of the path space $\{\gamma : [0,1] \to M \mid \gamma(0) \in L_{0},\gamma(1)\in L_{1}\}$.
\item Derived intersection theory \cite{ptvv}. Here $X$ is an $n$-shifted symplectic stack, and $f_{0}: L_{0} \to X$ and $f_{1} : L_{1} \to X$ are morphisms of stacks equipped with $n$-shifted Lagrangian structures, and we can form the derived fiber product
  \begin{equation*}
    \xymatrix{
      L_{0}\times_{X} L_{1} \ar[r]\ar[d] & L_{1} \ar[d] \\
      L_{0} \ar[r] & X
      }
    \end{equation*}
    which carries an $(n-1)$-shifted symplectic structure.
\end{enumerate}

In basic terms, our proposal composes these two ideas. Starting from an $(n-1)$-shifted symplectic stack $X$, with $n \geq 1$ (we do not consider negative shift), we iteratively apply (2), each time reducing the shift by one, until we obtain $0$-shifted Lagrangians in a $0$-shifted symplectic stack. Under the \emph{ansatz} that Lagrangian Floer cohomology can defined for those Lagrangians, we use (1) to define the ``top-level'' morphism spaces. If all goes well the result is an $(\infty,n)$-category that we denote $\nFuk{n}(X)$.

To make this more precise, we can present a list of axioms that implicitly and inductively characterize, at least partially, the object $\nFuk{n}(X)$, if it exists. We take $n \geq 1$ and $X$ an $(n-1)$-shifted symplectic stack.
\begin{enumerate}
\item (Base case) If $n = 1$, and $X$ is an ordinary symplectic manifold, then
  \begin{equation*}
    \nFuk{1}(X) = \Fuk(X),
  \end{equation*}
  where the latter is the Fukaya category defined using Floer theory, or an equivalent method.
\item (Induction step) If $n > 1$, the $(\infty,n)$-category $\nFuk{n}(X)$ has the following properties.
  \begin{enumerate}
  \item The objects of $\nFuk{n}(X)$ are $(n-1)$-shifted Lagrangians $L \to X$ (with extra data).
  \item Given two objects $L_{0}\to X$ and $L_{1}\to X$, the $(\infty,n-1)$-category of morphisms between them is given by
    \begin{equation*}
      \Hom_{\nFuk{n}(X)}(L_{0},L_{1}) = \nFuk{(n-1)}(L_{0}\times_{X}L_{1}),
    \end{equation*}
    where the fiber product is given the induced $(n-2)$-shifted symplectic structure.
  \item For three objects $L_{0},L_{1},L_{2}$, the binary composition of morphisms is given by a map
    \begin{equation*}
      \nFuk{(n-1)}(L_{0}\times_{X}L_{1}) \times \nFuk{(n-1)}(L_{1}\times_{X}L_{2}) \to \nFuk{(n-1)}(L_{0}\times_{X}L_{2})
    \end{equation*}
    that is induced by the $(n-1)$-shifted Lagrangian correspondence $L_{0}\times_{X} L_{1}\times_{X}L_{2}$. (The statement that this object is Lagrangian is a theorem of Ben-Bassat \cite{ben-bassat-lag-intersections}.)
  \end{enumerate}
\end{enumerate}

\begin{example}
  Our main examples have $n = 2$ and $X$ a $1$-shifted symplectic stack. In this case the objects of $\nFuk{2}(X)$ are $1$-shifted Lagrangians, and the Hom-category between $L_{0}\to X$ and $L_{1}\to X$ is the Fukaya category $\Fuk(L_{0}\times_{X} L_{1})$, where this is defined using Floer theory whenever possible.
\end{example}

Similar and closely related ideas have already appeared in the literature, including works of the Oxford group \cite{brav-et-al}, Bussi \cite{bussi-categorification}, Amorim--Ben-Bassat \cite{amorim-ben-bassat}, Calaque \cite{calaque-mapping-stacks},  Haugseng \cite{haugseng-spans}, Calaque-Haugseng-Scheimbauer \cite{calaque-haugseng-scheimbauer}, and Calaque-Safronov \cite{calaque-safronov}. The Oxford group usually considers taking $0$-shifted Lagrangian intersections in a $0$-shifted symplectic stack, and then adding an additional categorical layer by taking a perverse sheaf of vanishing cycles associated to the $(-1)$-shifted symplectic structure on the intersection. On the other hand, in the works of Calaque, Haugseng, Safronov, and Scheimbauer, the authors consider using iterated Lagrangian spans to build higher categories. This is something extremely close to our proposal, as it amounts to considering the  $(\infty,n-1)$-category defined by the same induction, but changing the base case so that the the \emph{set} of $1$-morphisms between two $0$-shifted symplectic stacks $X_{0}, X_{1}$ is the set of Lagrangian correspondences between them. Relative to this, our proposal adds an additional categorical layer by taking Floer cohomology of $0$-shifted Lagrangians.

\emph{Caveat lector.} This attempt at a characterization raises a number of questions, and seems fraught with technical difficulties.
\begin{enumerate}
\item The rigorous definition of this theory would seem to require the combination of Floer theory (or another theory of equivalent strength) with the derived geometry in the $C^{\infty}$ category. The work of Steffens \cite{steffens} is highly relevant in this regard, but there is clearly much more to be done.
\item The system of categories $\nFuk{n}(X)$, where $n \geq 1$ and $X$ is $(n-1)$-shifted symplectic, would constitute a vast generalization of the Weinstein $(\infty,2)$-category $\Symp$, for indeed we would have $\Symp \subseteq \nFuk{2}(\pt[1])$. The object $\Symp$ is currently under construction building on work of Wehrheim, Woodward, Bottman and others \cite{abouzaid-bottman}, and the construction of $\nFuk{n}(X)$ cannot be simpler than that of $\Symp$.
\item A more specific question that is likely to be easier to answer is, what are the ``extra data'' or ``brane structure'' that our $(n-1)$-shifted Lagrangians should carry in order to give rise to objects of $\nFuk{n}(X)$? The ``initial condition'' for this question is the definition of a brane structure in the $0$-shifted case, where a brane structure consists of a grading, a Pin structure, and in some formulations a local system on the Lagrangian. The higher brane structures on shifted Lagrangians must be deloopings of this notion, so that after taking derived intersection several times, we obtain the original notion of a $0$-shifted brane structure.
\end{enumerate}

Due to these difficulties, it seems reasonable to ask if the whole idea leads to anything interesting, since only then would the work to overcome the difficulties be justified. In the remainder of this paper, we consider several examples, and in each case try to indicate what may be gained by considering them in the light of the higher Fukaya category.

\section{The point}
Let $X = \pt$ be a point. The tangent complex $\bT_{\pt}$ and the cotangent complex $\bL_{\pt}$ are both zero, and so the zero map is a quasi-isomorphism $\bT_{\pt} \to \bL_{\pt}[n]$ for any $n$. Hence the zero $n$-shifted $2$-form is symplectic. We denote by $\pt[n]$ the point considered as an $n$-shifted symplectic manifold this way.

An important theorem of \cite{ptvv} states that Lagrangian structures on the terminal morphism $f : L \to \pt[n]$ are equivalent to $(n-1)$-shifted symplectic structures on $L$.

In the case $n=1$, we see that any ordinary symplectic manifold $(M,\omega)$ defines a $1$-shifted Lagrangian in $\pt[1]$, and
\begin{equation*}
  \Hom_{\nFuk{2}(\pt[1])}(M_{0},M_{1}) = \Fuk(M_{0} \times_{\pt[1]} M_{1}) = \Fuk(\overline{M}_{0} \times M_{1}) = \Hom_{\Symp}(M_{0},M_{1})
\end{equation*}
where the bar denotes the reversal of the symplectic form. This yields the aforementioned embedding $\Symp \subseteq \nFuk{2}(\pt[1])$. The larger category contains many $0$-shifted symplectic stacks that are not smooth manifolds, and we do not know how to define the Fukaya category for such objects in general. Still, given the rapid development of derived geometry, it does not seem unlikely that such a theory will eventually exist.\footnote{The work of Lekili and Segal \cite{lekili-segal} can be viewed as defining the Fukaya category of certain non-smooth $0$-shifted symplectic stacks, namely the reductions by a torus action at critical level. See Section \ref{sec:coadjoint} for further discussion of this example.}

For higher $n$, we can at least make a few formal remarks. The category $\nFuk{(n+1)}(\pt[n])$ is symmetric monoidal under Cartesian product, and the unit is $\pt[n-1]$. The endomorphisms of the unit object are
\begin{equation*}
  \Hom_{\nFuk{(n+1)}(\pt[n])}(\pt[n-1],\pt[n-1]) = \nFuk{n}(\pt[n-1]).
\end{equation*}
Any $(n-1)$-shifted symplectic stack $X$ is an object of $\nFuk{(n+1)}(\pt[n])$, and we have
\begin{equation*}
  \Hom_{\nFuk{(n+1)}(\pt[n])}(\pt[n-1],X) = \nFuk{n}(X).
\end{equation*}

If we denote $\nSymp{n} = \nFuk{(n+2)}(\pt[n+1])$  (so that $\nSymp{0} \supseteq \Symp$), these remarks say, firstly, that loops on $\nSymp{n}$ is equivalent to $\nSymp{(n-1)}$, and secondly, that $\nSymp{n}$ contains all of the structure of the higher Fukaya categories of all shifted symplectic stacks with shift between $0$ and $n$.

\section{The one-shifted cotangent stack}
Let $M$ be a smooth manifold, and let $n \geq 0$. The $n$-shifted cotangent stack is defined by shifting the degree of the fibers of the cotangent bundle:
\begin{equation*}
  \rT^{*}[n]M = \Spec_{M}(\Sym_{\cO_{M}}(\bT_{M}[-n])),
\end{equation*}
where this formula should be interpreted in terms of $C^{\infty}$-derived geometry: $\bT_{M}$ is the tangent sheaf, $\cO_{M}$ is structure sheaf of $C^{\infty}$ functions, and $\Spec_{M}$ is the relative spectrum of a sheaf of $C^{\infty}$-rings. The stack $\rT^{*}[n]M$ is $n$-shifted symplectic; in the algebraic category this is a theorem of Calaque \cite{calaque-cotangent}.

We consider the case $n=1$ for the remainder of this section. We shall investigate some natural objects in $\nFuk{2}(\rT^{*}[1]M)$.

\subsection{The zero section}

The zero section $s : M \to \rT^{*}[1]M$ carries a canonical $1$-shifted Lagrangian structure; again, in the algebraic category this is a theorem of Calaque \cite{calaque-cotangent}. The derived self-intersection of this object is
\begin{equation*}
  M \times_{\rT^{*}[1]M}M \cong \rT^{*}M,
\end{equation*}
the ordinary cotangent bundle of $M$, equipped with the canonical symplectic form $\omega_{\mathrm{can}}$. Thus the category of endomorphisms of the zero section as an object of $\nFuk{2}(\rT^{*}[1]M)$ is
\begin{equation*}
  \Hom_{\nFuk{2}(\rT^{*}[1]M)}(M,M) = \Fuk(\rT^{*}M).
\end{equation*}
This category must be monoidal, since it is the endomorphism category of an object. Indeed, there is a natural Lagrangian correspondence in $\overline{\rT^{*}M} \times \overline{\rT^{*}M} \times \rT^{*}M$ that is the graph of fiberwise addition in the cotangent bundle, which is expected to define a monoidal structure on $\Fuk(\rT^{*}M)$, as first considered by Subotic \cite{subotic}, see also \cite{pascaleff,abouzaid-bottman,abouzaid-bottman-niu}.

We can also obtain other objects by modifying the structure of the Lagrangian morphism $M \to \rT^{*}[1]M$. According to \cite{ptvv}, and adopting their notation, starting with the canonical Lagrangian structure on the morphism $M\to \rT^{*}[1]M$, we can modify it by adding a class in $\pi_{1}(\cA^{2,\mathrm{cl}}(M,1))$. Since $M$ is a smooth manifold, this group boils down to $\Omega^{2,\mathrm{cl}}(M)$, the space of (ordinary) closed $2$-forms on $M$. Generally speaking, such a modification is merely an isotropic structure, but in the case at hand it is always Lagrangian (i.e., the nondegeneracy condition is satisfied). Given $\beta \in \Omega^{2,\mathrm{cl}}(M)$ we denote by $M_{\beta}$ the Lagrangian structure on $M \to \rT^{*}[1]M$ obtained by this modification.

We then have
\begin{equation*}
  M \times_{\rT^{*}[1]M}M_{\beta} \cong (\rT^{*}M,\omega_{\mathrm{can}}+\pi^{*}\beta),
\end{equation*}
where the right-hand side is the \emph{twisted} or \emph{magnetic} cotangent bundle. Hence the category of morphisms from $M$ to $M_{\beta}$ is $\Fuk(\rT^{*}M,\omega_{\mathrm{can}}+\pi^{*}\beta)$. This category is therefore a module category over the monoidal category $\Fuk(\rT^{*}M)$, where the action is again given by the fiberwise addition correspondence. (More generally, the fiberwise addition correspondence adds the twists.) The symplectic cohomology of these twisted cotangent bundles has been studied \cite{benedetti-ritter,groman-merry}

\begin{remark}
  One might also attempt to obtain more Lagrangians by deforming the morphism $M \to \rT^{*}[1]M$ away from the zero section. In the present setting this is not possible, for the following reason. Sections $M \to \rT^{*}[1]M$ are precisely $1$-shifted $1$-forms on $M$. Since $M$ is classical (not stacky), \cite[Proposition 1.14]{ptvv} implies that such objects are elements of $\rH^{1}(M,\Omega_{M}^{1})$. Since $M$ is a smooth manifold, it is affine as a $C^{\infty}$-scheme, and hence this group vanishes.
\end{remark}

\subsection{Symplectic fibrations}

Let $E$ be a smooth manifold and let $\pi : E \to M$ be a smooth fibration (i.e., a submersion). Composing with the zero section $s$ gives a morphism $s\pi : E \to \rT^{*}[1]M$, and we may ask what are Lagrangian structures on this morphism.

It turns out that such Lagrangian structures are related to closed $2$-forms on $E$ that are symplectic on the fibers, that is, the set of forms
\begin{equation*}
  \{\tau \in \Omega^{2}(E) \mid \text{$d\tau = 0$ and $\tau|_{\pi^{-1}(x)}$ is nondegenerate}\}.
\end{equation*}
Such a form $\tau$ makes all fibers $\pi^{-1}(x)$ symplectic, and it also induces a symplectic connection on $\pi : E \to M$. This means that all fibers are pairwise symplectomorphic, and so the structure group is reduced to $\Symp(\pi^{-1}(x))$. If $M$ is simply connected, the structure group reduces further to $\Ham(\pi^{-1}(x))$. The technical name for this structure is a \emph{locally Hamiltonian fibration}, but we shall also call it a \emph{symplectic fibration}.

\begin{theorem}
  $1$-shifted Lagrangian structures on the morphism $s\pi : E \to \rT^{*}[1]M$ are equivalent to symplectic fibration structures on $\pi : E \to M$, that is, elements $\tau \in \Omega^{2,\mathrm{cl}}(E)$ such that $\tau$ is nondegenerate when restricted to the fibers of $\pi$.
\end{theorem}

\noindent \emph{Sketch of proof.} A Lagrangian structure implies a duality between the tangent complex $\bT_{E}$ and a shift of the relative tangent complex $\bT_{s\pi}$ of the morphism, so the first step is to investigate these complexes. 

The tangent complex of $E$ sits in a triangle
\begin{equation*}
  \bT_{\pi} \to \bT_{E} \to \pi^{*}\bT_{M}.
\end{equation*}
Since $\pi : E \to M$ is a submersion of smooth manifolds, this is just a short exact sequence of vector bundles on $E$. Note that $\bT_{\pi}$ is the vertical tangent bundle $\ker(D\pi)$, whose fibers are the tangent spaces to the fibers of $\pi$.

The relative tangent complex $\bT_{s\pi}$ sits a triangle
\begin{equation*}
  \bT_{\pi} \to \bT_{s\pi} \to \pi^{*}\bT_{s}
\end{equation*}
coming from the fact that $s\pi$ is a composition of two maps. The relative tangent complex $\bT_{s}$ sits in a triangle
\begin{equation*}
  \bT_{s} \to \bT_{M} \to s^{*}\bT_{\rT^{*}[1]M}.
\end{equation*}
Because we are on the zero section, there is a canonical splitting $s^{*}\bT_{\rT^{*}[1]M} \cong \bT_{M} \oplus \bL_{M}[1]$ as the sum of the tangent and shifted cotangent bundles. The map from $\bT_{M}$ to this object is the inclusion of the first summand, so we deduce that the triangle
\begin{equation*}
  \bT_{s} \to 0 \to \bL_{M}[1]
\end{equation*}
is exact, giving an isomorphism $\bT_{s} \cong \bL_{M}$.

The outcome of the previous two paragraphs is that the tangent and relative tangent complexes have a similar structure:
\begin{align*}
    \bT_{\pi} \to & \bT_{E} \to \pi^{*}\bT_{M},\\
    \bT_{\pi} \to & \bT_{s\pi} \to \pi^{*}\bL_{M}.
  \end{align*}
A Lagrangian structure must induce a duality between $\bT_{E}$ and $\bT_{s\pi}$. We observe that, while the terms $\pi^{*}\bT_{M}$ and $\pi^{*}\bL_{M}$ are naturally dual, we would seem to need a pairing on $\bT_{\pi}$. Also, we need a splitting of one of the exact sequences so as to match the direction of one to the dual of the other. Fortunately, a symplectic fibration form $\tau$ provides exactly these data.

Starting from a form $\tau \in \Omega^{2,\mathrm{cl}}(E)$ such that $\tau$ is nondegenerate on $\bT_{\pi}$, we construct the candidate Lagrangian structure as follows. The morphism $s : M \to \rT^{*}[1]M$ has a canonical Lagrangian structure, and pulling this back under $\pi : E \to M$ we obtain at least an isotropic structure on $s\pi$. Outside of trivial cases, this isotropic structure is never Lagrangian, because it is degenerate along $\bT_{\pi}$. However, we can modify the isotropic structure by adding $\tau \in \Omega^{2,\mathrm{cl}}(E) = \pi_{1}(\cA^{2,\mathrm{cl}}(E,1))$, where we are again using the fact that $E$ is a smooth manifold. As long as $\tau$ is nondegenerate on $\bT_{\pi}$, it determines a $\tau$-orthogonal splitting $\bT_{E} \cong \bT_{\pi} \oplus \pi^{*}\bT_{M}$. This fixes both of the problems in the previous paragraph, so $\tau$ defines a nondegenerate pairing between $\bT_{E}$ and $\bT_{s\pi}$. Thus we obtain a Lagrangian structure. This ends the sketch of the proof. 

\begin{remark}
  When $M = \pt$, this theorem reduces to the statement that Lagrangian structures on $E \to \pt[1]$ are symplectic structures on $E$ \cite{ptvv}.
\end{remark}

\begin{example}
Let us now consider $(\pi : E \to M, \tau)$ as an object of $\nFuk{2}(\rT^{*}[1]M)$. Taking the morphisms from the zero section yields
\begin{equation*}
 \Hom_{\nFuk{2}(\rT^{*}[1]M}(M, (\pi : E \to M, \tau)) = \Fuk(M \times_{\rT^{*}[1]M} E).
\end{equation*}
It turns out that the object $Z = M \times_{\rT^{*}[1]M} E$ is a smooth symplectic manifold. The smooth structure is given by writing $Z$ as the fiber product 
\begin{equation*}
  \xymatrix{
    Z \ar[r]\ar[d] & \rT^{*}M \ar[d] \\
    E \ar[r] & M
  }
\end{equation*}
and the symplectic form is the sum of the pullbacks of $\tau$ and the canonical symplectic form on $\rT^{*}M$. (When $M = S^{1}$, the object $Z$ is known as a symplectic mapping torus, the Fukaya category of which has been studied by Kartal \cite{kartal-dynamical,kartal-distinguishing}.) Once again, the $2$-categorical structure implies that the category $\Fuk(Z)$ is a module category for the monoidal category $\Fuk(\rT^{*}M)$.
\end{example}
We now wish to connect to the notion of \emph{local systems of categories} over $M$. These objects are studied in great detail  in recent joint work of Emanuele Pavia and the present authors  \cite{pps-loc-sys}.  The most direct definition of a local system of categories involves considering $M$ as an $\infty$-groupoid, which is a type of $\infty$-category, and then taking functors from this $\infty$-groupoid to $\dgcat$, the $\infty$-category of $k$-linear differential graded categories. The $(\infty,2)$-category of local systems of categories on M is then
\begin{equation*}
  \locsyscat(M) = \Fun(M,\dgcat).
\end{equation*}
Alternatively, we can take the Betti stack of $M$, denoted $M_B$, and take quasicoherent sheaves of categories in the sense of Gaitsgory \cite{gaitsgory},
\begin{equation*}
  \locsyscat(M) \cong \shvcat(M_{B}).
\end{equation*}
We now apply some of the concepts from \cite{gaitsgory}. Let $\underline{\Perf(k)}$ denote the trivial local system of categories on $M$; its endomorphism category is
\begin{equation*}
  \Hom_{\locsyscat(M)}(\underline{\Perf(k)},\underline{\Perf(k)}) \cong \locsys(M),
\end{equation*}
the $(\infty,1)$-category of local systems of $k$-vector spaces on $M$. Notably $\locsys(M)$ is equivalent to the wrapped Fukaya category of $\rT^{*}M$, by a theorem of Abouzaid \cite{abouzaid-fiber-generates,abouzaid-based-loops}.

Furthermore, given any local system of categories $\cE$, the category of global sections is a module category for the endomorphism category of $\underline{\Perf(k)}$:
\begin{equation*}
  \Gamma(M,\cE) = \Hom_{\locsyscat(M)}(\underline{\Perf(k)},\cE) \in \Mod_{\locsys(M)}.
\end{equation*}
This actually arises from a functor
\begin{equation*}
  \Gamma^{\mathrm{enh}} : \locsyscat(M) \to \Mod_{\locsys(M)}
\end{equation*}
This functor has a left adjoint, that we shall denote $L : \Mod_{\locsys(M)} \to \locsyscat(M)$.

\begin{remark}
  The adjunction
  \begin{equation*}
    L : \Mod_{\locsys(M)}   \rightleftarrows \locsyscat(M) : \Gamma^{\mathrm{enh}}
  \end{equation*}
  is an adjoint equivalence precisely when the Betti stack $M_{B}$ is \emph{$1$-affine}. This condition rarely holds for Betti stacks, see \cite{pps-koszul} for a detailed analysis of this problem.
\end{remark}

Let us draw out the connection between symplectic fibrations and local systems of categories that has already been hinted at. For this it is essential that we use wrapped Fukaya categories throughout. We propose the following:
\begin{enumerate}
\item A symplectic fibration $(\pi : E \to M,\tau)$ corresponds to a local system of categories $\cE$ over $M$.
\item The trivial fibration $(M \to M, \tau = 0)$, or equivalently the zero section $s : M \to \rT^{*}[1]M$, corresponds to the trivial local system $\underline{\Perf(k)}$.
\item By taking endomorphisms, there is a monoidal equivalences between $\Fuk(\rT^{*}M)$ and $\locsys(M)$.
\item The category $\Fuk(M \times_{\rT^{*}[1]M} E)$ is equivalent to $\Gamma(M,\cE)$, in a way compatible with the module structures.
\end{enumerate}
Thus we expect there to be a diagram as follows.
\begin{equation*}
  \xymatrix{
    \nFuk{2}(\rT^{*}[1]M) \ar[r] & \Mod_{\locsys(M)} \ar[d]_{L}\\
    & \locsyscat(M) \ar@<-1ex>[u]_{\Gamma^{\mathrm{enh}}}
  }
\end{equation*}
Where the top arrow is taking Hom with the zero section, and the vertical arrows are Gaitsgory's adjunction. It would be nice to be able fill in the diagram with a diagonal arrow. Taking a clue from the theory of constructible sheaves, we would expect $\locsyscat(M)$ to arise from shifted Lagrangians that are ``supported on the zero section,'' as $s\pi : E \to \rT^{*}[1]M$ is. This motivates a conjecture.
\begin{conjecture}
  There is a fully faithful embedding of $(\infty,2)$-categories
  \begin{equation*}
    \locsyscat(M) \to \nFuk{2}(\rT^{*}[1]M).
  \end{equation*}
  The essential image contains the objects arising from symplectic fibrations, and it is contained in a subcategory $\nFuk{2}_{M}(\rT^{*}[1]M)$ whose objects are $1$-shifted Lagrangians that factor through the zero section.
\end{conjecture}

\subsection{Lefschetz fibrations}

We now generalize the previous example to allow the map $\pi : E \to M$ to have singularities, for instance we could consider the case where $\pi : E \to \bC$ is a Lefschetz fibration.

If the map $\pi : E \to M$ is not a submersion, then the morphism $s\pi : E \to \rT^{*}[1]M$ does not admit a Lagrangian structure. This is true even at the formal level, in the sense that the complexes $\bL_{E}$ and $\bT_{s\pi}$ are not isomorphic (a straightforward calculation shows that, at the critical points, the stalks of their cohomology sheaves have different ranks). There is a fix for this problem that is, as one says, obvious once you see it.

The \emph{singularity stack} of the map $\pi$ is the stack
\begin{equation*}
  \Sing(\pi) = \Spec_{E}(\Sym_{\cO_{E}}(\rH^{1}(\bT_{\pi})[-1]))
\end{equation*}
The space $\rH^{1}(\bT_{\pi})$ is the cokernel of $D\pi$, and this space lives in cohomological degree 1, which is its natural cohomological degree. At all points where $\pi$ is submersive, we have $\rH^{1}(\bT_{\pi}) = 0$, so $\Sing(\pi)$ differs from $E$ only at the critical points. If we stratify $E$ according to the rank of $D\pi$, then we could view $\Sing(\pi)$ as a stratified graded vector bundle.

\begin{remark}
  The definition of $\Sing(\pi)$ is similar in spirit to definition of the singular support variety of Arinkin-Gaitsgory \cite{arinkin-gaitsgory}, but with a different shift.  
\end{remark}

The singularity stack has a natural morphism to the one-shifted cotangent stack, constructed as follows. Recall the definition
\begin{equation*}
  \rT^{*}[1]M = \Spec_{M}(\Sym_{\cO_{M}}(\bT_{M}[-1])).
\end{equation*}
There is a natural map of complexes $\pi^{*}\bT_{M}[-1] \to \rH^{1}(\bT_{\pi})[-1]$, coming from the definition of $\rH^{1}(\bT_{\pi})$ as the cokernel of a map toward $\bT_{M}$. Taking relative $\Spec$ gives a map in the opposite direction,
\begin{equation*}
  \widetilde{\pi} : \Sing(\pi) \to \rT^{*}[1]M.
\end{equation*}
We propose to put $1$-shifted Lagrangian structures on this morphism, so as to obtain objects of $\nFuk{2}(\rT^{*}[1]M)$.

We now recall one notion of a symplectic fibration with singularities \cite{seidel-fcplt}. For this we suppose that $\pi : E \to M$ has only isolated critical points, and that, the critical points are fully degenerate, so $D\pi_{p} = 0$, when $p$ is critical. (Lefschetz fibrations are an example.) A \emph{symplectic fibration with singularities} structure is induced by a closed two-form $\tau\in \Omega^{2,\mathrm{cl}}(E)$ such that $\tau$ is symplectic when restricted to the smooth loci of the fibers of $\pi$, and moreover, $\tau$ is symplectic in a neighborhood of each critical point.

More generally, we could consider closed $2$-forms  whose restriction the the singular (non-constant-rank) distribution $\ker(D\pi)$ is symplectic.

\begin{theorem}
  \label{thm:sing}
  Suppose that $\pi : E \to M$ is a fibration with singularities. If $\tau \in \Omega^{2,\mathrm{cl}}(E)$ is a closed $2$-form  whose restriction to the singular distribution $\ker(D\pi)$ is symplectic, then $\tau$ induces a Lagrangian structure on the morphism $\widetilde{\pi} : \Sing(\pi) \to \rT^{*}[1]M$.
\end{theorem}

\noindent \emph{Sketch of proof.} At points of $E$ where $D\pi$ has full rank, the analysis of the tangent and relative tangent complexes is the same as before, so we focus on a critical point $p \in E$, with critical value $v = \pi(p) \in M$.

For the tangent complex of the source, we have
\begin{equation*}
  \bT_{\Sing(\pi),p} \cong \bT_{E,p} \oplus \rH^{-1}(\bL_{\pi,p})[1] 
\end{equation*}
where the second term lives in the shifted cotangent direction, and there is no differential in this presentation. For later use we shall rewrite this as
\begin{equation*}
  \bT_{\Sing(\pi),p} = \bT_{E,p} \oplus \ker(D\pi^{*})[1] = \ker(D\pi) \oplus C \oplus \ker(D\pi^{*})[1]
\end{equation*}
where $C \subseteq \bT_{E,p}$ is a subspace complementary to the kernel of $D\pi$. We then have $C \cong \im(D\pi)$.

The relative tangent complex of $s\pi$ is given by
\begin{equation*}
  \bT_{s\pi,p} \cong \bT_{\pi,p} \oplus \pi^{*}\bL_{M,v}
\end{equation*}
In passing from $s \pi$ to $\widetilde{\pi}$, the tangent complex to the source obtains a new term $\rH^{-1}(\bL_{\pi,p})[1]$, and there is a differential that includes this factor into the $\pi^{*}\bL_{M,v}$. Thus as a complex we have 
\begin{equation*}
  \bT_{\widetilde{\pi},p} \cong\{ \rH^{-1}(\bL_{\pi,p}) \to \bT_{E,p} \oplus \pi^{*}\bL_{M,v} \to \bT_{M,v}\},
\end{equation*}
where the differential connects $\bT_{E,p}$ to $\bT_{M,v}$ via $D\pi$, and $\rH^{-1}(\bL_{\pi,p})$ to $\pi^{*}\bL_{M,v}$ via the natural inclusion. The cohomology sheaves of $\bT_{\widetilde{\pi},p}$ are then
\begin{align*}
  \rH^{0}(\bT_{\widetilde{\pi},p}) &\cong \ker(D\pi) \oplus (\pi^{*}\bL_{M,v}/\rH^{-1}(\bL_{\pi,p}))\\
  \rH^{1}(\bT_{\widetilde{\pi},p}) &\cong \rH^{1}(\bT_{\pi,p}) = \coker(D\pi). 
\end{align*}
The space $\pi^{*}\bL_{M,v}/\rH^{-1}(\bL_{\pi,p})$ is isomorphic to the image of $D\pi^{*}$ by the first isomorphism theorem. Thus we have
\begin{equation*}
  \bT_{\widetilde{\pi},p} \cong \ker(D\pi) \oplus \im(D\pi^{*}) \oplus \coker(D\pi)[-1].
\end{equation*}
Comparing $\bT_{\Sing(\pi),p}$ with $\bT_{\widetilde{\pi},p}$, we see that the terms $\ker(D\pi^{*})[1]$ and $\coker(D\pi)[-1]$ are canonically dual, while the ``horizontal'' subspace $C$ pairs perfectly with $\im(D\pi^{*})$. So it suffices that the $2$-form $\tau$ induces a nondegenerate pairing on the remaining terms, which are both  $\ker(D\pi)$. This ends the sketch of the proof.

\begin{example}
Suppose $\pi : E \to \bC$ is a Lefschetz fibration with a single critical point $p$ such that $\pi(p) = 0$. Then
\begin{equation*}
  \Sing(\pi) = E \cup_{p} \rT^{*}_{0}[1]\bC,
\end{equation*}
and the map $\widetilde{\pi}$ maps $E$ to $\bC \subset \rT^{*}[1]\bC$ and the cotangent fiber to itself.
\end{example}

\begin{example}
More generally, suppose $\pi: E \to M$ has isolated critical points with distinct critical values. Then there is a factorization of $\widetilde{\pi}$
\begin{equation*}
  \Sing(\pi) \to M \cup \bigcup_{p \in \crit(\pi)} \rT^{*}_{\pi(p)}[1]M \to \rT^{*}[1]M,
\end{equation*}
where the second map is an embedding of a conical Lagrangian substack that is the union of the zero-section and a set of shifted cotangent fibers \cite{calaque-cotangent}. (This phenomenon, where one Lagrangian morphism factors through another, is an interesting feature of derived symplectic geometry.)

This picture strongly suggests an analogy to constructible sheaves and singular support. Just as a complex of constructible sheaves on $M$ has a singular support that is a conical Lagrangian subvariety of $\rT^{*}M$, a $1$-shifted Lagrangian coming from a symplectic fibration with singularities has a ``singular support'' that is a conical $1$-shifted Lagrangian in $\rT^{*}[1]M$. In the case where $\pi : E \to M$ is a smooth fibration, this singular support reduces to the zero section $M$, analogous to a local system.


\end{example}

The fundamental work of Seidel \cite{seidel-fcplt} shows that, given a symplectic Lefschetz fibration $\pi : E \to S$, where $S$ is a Riemann surface, we obtain (what has subsequently come to be known as) a \emph{schober} on $S$, which has singularities at $\critv(\pi)$, the set of critical values of $\pi$. A schober is a the categorification of the concept of a perverse sheaf, and so it generalizes local systems of categories. On the other hand, we have just seen that one obtains a $1$-shifted Lagrangian in $\rT^{*}[1]S$ from these same data, which is supported on the $1$-shifted conormal bundle $\rT^{*}_{\critv(\pi)}[1]S$ of the set of critical values. 
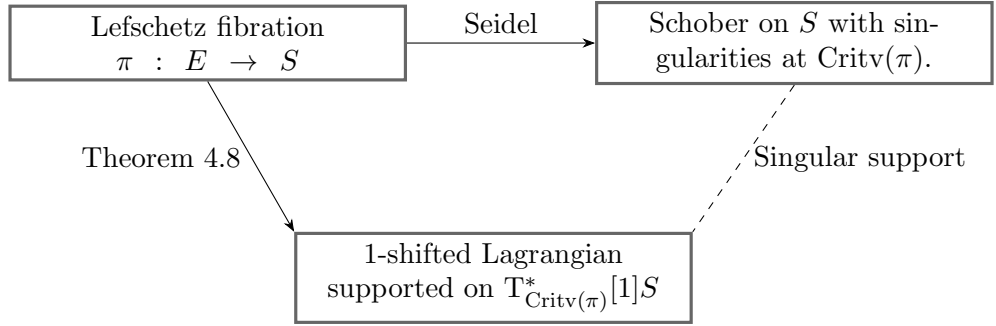
\begin{figure}[!h]
  \centering
  \begin{tikzpicture}
    [
    squarednode/.style={%
      rectangle,
      draw=black!60,
      fill=white,
      very thick,
      minimum size=5mm,
      text centered,
      text width=5cm,
      node distance=2.5cm
    }
    ]
    \node[squarednode]      (lefschetz)                             {Lefschetz fibration\\ $\pi : E \to S$};
    \node[squarednode]      (schober)  [right=of lefschetz] {Schober on $S$ with singularities at $\critv(\pi)$. };
    \node[squarednode]      (lag)  [below right=2cm and -1.5cm of lefschetz]   {$1$-shifted Lagrangian supported on $\rT^{*}_{\critv(\pi)}[1]S$};
    
    \draw[-Stealth] (lefschetz.east) -- node [above,midway] {Seidel} (schober.west);
    \draw[-Stealth] (lefschetz.south) -- node [left,midway] {Theorem \ref{thm:sing}}(lag.north west);
    \draw[dashed] (schober.south) -- node [right,midway] {Singular support}(lag.north east);
  \end{tikzpicture}
  \caption{Structures deriving from a Lefschetz fibration.}
  \label{figure:link}
\end{figure}

This suggests that a schober on $S$ should be directly related to a $1$-shifted Lagrangian in $\rT^{*}[1]S$. This would of course be a generalization, at the categorified level, of the notion of singular support for a perverse sheaf. Figure \ref{figure:link} displays these relationships.

One can then ask whether the relationship between schobers and the $1$-shifted Lagrangians holds at the $2$-categorical level. This leads to a conjectural categorification of the Nadler-Zaslow correspondence.
\begin{conjecture}
  Let $\cS$ be a stratification of $M$. Then there is a fully faithful embedding of $(\infty,2)$-categories
  \begin{equation*}
    \schobers(M,\cS) \to \nFuk{2}(\rT^{*}[1]M)
  \end{equation*}
  whose essential image is contained in the subcategory of objects supported on $\rT^{*}_{\cS}[1]M$, the shifted conormal to the stratification.
\end{conjecture}
 
\section{The coadjoint stack}
\label{sec:coadjoint}

Let $G$ be a compact Lie group, and consider the stack
\begin{equation*}
  X = [\frakg^{*}/G] = \rT^{*}[1]\rB G,
\end{equation*}
which is the quotient stack of the coadjoint action of $G$ on $\frakg^{*}$, and also the $1$-shifted cotangent bundle of $\rB G = [\pt/G]$. The latter presentation makes it clear that $[\frakg^{*}/G]$ is $1$-shifted symplectic.

\subsection{Hamiltonian $G$-manifolds}
\label{sec:ham-g-man}

The shifted symplectic geometry of $[\frakg^{*}/G]$ recovers many important results in the theory of Hamiltonian $G$-manifolds \cite{safronov-quasi-ham,safronov-geometric-quant}. So in a certain sense, this theory is already well-studied. Suppose $Y$ is a symplectic manifold with a Hamiltonian $G$-action, and moment map $\mu: Y \to \frakg^{*}$. Then $\mu$ induces a morphism
\begin{equation*}
  [Y/G] \to [\frakg^{*}/G]
\end{equation*}
that naturally admits a $1$-shifted Lagrangian structure. The reduced space is obtained by derived fiber product \cite{safronov-quasi-ham}
\begin{equation}
  \label{eq:reduction}
  Y/\!\!/_{0}G \cong [\pt/G] \times_{[\frakg^{*}/G]} [Y/G].
\end{equation}
This can also be viewed as the derived intersection of the zero section in $\rT^{*}[1]\rB G$ with the Lagrangian $[Y/G]$.
To be clear, the isomorphism \eqref{eq:reduction} holds when $0$ is a regular value of $\mu$, making the left-hand side smooth. When $0$ is a critical value, we can use \eqref{eq:reduction} as the \emph{definition} of the symplectic reduction, and in that case it is a non-smooth $0$-shifted symplectic stack.

\begin{remark}
 In the case where $G$ is a torus, the Fukaya categories of symplectic reductions at critical values of the moment map have been studied by Lekili-Segal \cite{lekili-segal}. Their work thus provides a beginning indication of what the Fukaya category of a $0$-shifted symplectic stack may look like.
\end{remark}

The cotangent fiber of $\rT^{*}[1]\rB G$ corresponds to the quotient morphism $\frakg^{*} \to [\frakg^{*}/G]$. We have the derived intersections
\begin{align*}
  \frakg^{*} \times_{[\frakg^{*}/G]} \frakg^{*} &\cong \rT^{*}G,\\
  \frakg^{*} \times_{[\frakg^{*}/G]} [Y/G] &\cong Y.
\end{align*}
Formally it follows that
\begin{align*}
  \Hom_{\nFuk{2}([\frakg^{*}/G])}(\frakg^{*}, \frakg^{*}) &\cong \Fuk(\rT^{*}G),\\
  \Hom_{\nFuk{2}([\frakg^{*}/G])}(\frakg^{*}, [Y/G]) &\cong \Fuk(Y),
\end{align*}
implying that $\Fuk(\rT^{*}G)$ is monoidal, and that $\Fuk(Y)$ is a module category. In some cases, this expectation has been rigorously justified, and moreover applied \cite{evans-lekili}. 

\subsection{Teleman's action}
\label{teleact}
We shall not attempt here to rephrase every result about Hamiltonian $G$-manifolds in the $1$-shifted framework, but there is an important point where the higher Fukaya category concept seems to be clarifying. A remarkable result of Teleman \cite{teleman-icm,teleman-coulomb} states that the $G$-equivariant quantum cohomology $\QH_{G}^{*}(Y)$ is a module over the ring $\rH_{*}^{G}(\Omega G)$, where this ring plays a fascinating role in geometric representation theory as $\cO(\BFM(G^{\vee}))$, the ring of functions on the Bezrukavnikov-Finkelberg-Mirkovic space of the Langlands dual group. This action has been rigorously constructed by González-Mak-Pomerleano \cite{gonzalez-mak-pomerleano-nil,gonzalez-mak-pomerleano-coulomb}, and it refines both the Seidel representation \cite{seidel-rep} and its higher-degree versions constructed by Savelyev \cite{savelyev-hofer}.

Setting aside the question of rigor, which we have no claim to, we shall argue that the theory of higher Fukaya categories gives a reasonably convincing heuristic for the existence of an action of $\rH_{*}^{G}(\Omega G)$ on $\QH_{G}^{*}(Y)$, though we need to make a number of (possibly optimistic) assumptions. Thus this notable result can be naturally situated in the theory of higher Fukaya categories as we understand them. We view this fact as favoring our claim that our version of this category is indeed relevant to 3D HMS and related matters.

The calculations we shall suggest live in the $3$-category that contains $[\frakg^{*}/G]$ as an object,
\begin{equation*}
  [\frakg^{*}/G] \in \nFuk{3}(\pt[2]) = \nSymp{1}.
\end{equation*}
We present the argument as follows.
\begin{enumerate}
\item The diagonal $\Delta \subset [\frakg^{*}/G] \times [\frakg^{*}/G]$ corresponds to the identity morphism of $[\frakg^{*}/G]$ in $\nSymp{1}$.
\item The Hamiltonian $G$-space $[Y/G] \to [\frakg^{*}/G]$ is a $1$-morphism from $\pt[1]$ to $[\frakg^{*}/G]$; we denote this morphism $L_{Y}$.
\item The equation $\Delta \circ L_{Y} = L_{Y}$ holds as a composition of $1$-morphisms. By taking $1$-categories of $2$-endomorphisms, we thus obtain an action
  \begin{equation*}
    \nHom{2}(\Delta,\Delta) \times \nHom{2}(L_{Y},L_{Y}) \to \nHom{2}(L_{Y},L_{Y}).
  \end{equation*}
\item The $1$-categories of $2$-endomorphisms in the preceding point are themselves monoidal (under composition of $2$-morphisms), and the action map is a morphism of monoidal categories (by the exchange relation). Hence the action maps preserve the monoidal unit objects.
\item By taking the $3$-endomorphisms of the monoidal unit objects, we obtain a action
  \begin{equation*}
    \nHom{3}(\bone_{\Delta},\bone_{\Delta}) \times \nHom{3}(\bone_{L_{Y}},\bone_{L_{Y}}) \to \nHom{3}(\bone_{L_{Y}},\bone_{L_{Y}}).
  \end{equation*}
  We have now used all levels of the $3$-categorical structure. It remains to determine the $E_{3}$-algebra $\nHom{3}(\bone_{\Delta},\bone_{\Delta})$ and the $E_{2}$-algebra $\nHom{3}(\bone_{L_{Y}},\bone_{L_{Y}})$.
\item To understand the $E_{2}$-algebra $\nHom{3}(\bone_{L_{Y}},\bone_{L_{Y}})$, we go back to the monoidal $1$-category
  \begin{equation*}
    \nHom{2}(L_{Y},L_{Y}) = \nFuk{1}([Y/G] \times_{[\frakg^{*}/G]} [Y/G]).
  \end{equation*}
  The right-hand side is the Fukaya category of the $0$-shifted symplectic stack; since this stack is only rarely a classical smooth manifold, we find ourselves outside the domain where Floer cohomology is usually defined. Nevertheless, we intend to take the Lagrangian Floer cohomology of the diagonal with itself in $[Y/G] \times_{[\frakg^{*}/G]} [Y/G]$. Our heuristic for this is based on the fact that, when working nonequivariantly, the Lagrangian Floer cohomology of the diagonal in $Y \times Y$ is the quantum cohomology $\QH^{*}(Y)$ (we are assuming $Y$ is compact; otherwise we would consider symplectic cohomology). Thus we expect
  \begin{equation*}
    \nHom{3}(\bone_{L_{Y}},\bone_{L_{Y}}) = \HF^{*}_{[Y/G] \times_{[\frakg^{*}/G]} [Y/G]}(\Delta,\Delta) = \QH^{*}_{G}(Y).
  \end{equation*}
  
\item To understand the $E_{3}$-algebra $\nHom{3}(\bone_{\Delta},\bone_{\Delta})$, observe that the derived self-intersection of the diagonal $\Delta \subset [\frakg^{*}/G] \times [\frakg^{*}/G]$ is the stack of (locally constant) loops
  \begin{equation*}
    L[\frakg^{*}/G] = L\rT^{*}[1]\rB G = \rT^{*}(L\rB G) = \rT^{*}[G/G],
  \end{equation*}
  which we identify as the cotangent bundle of the group adjoint quotient stack. Thus we have
  \begin{equation*}
    \nHom{2}(\Delta,\Delta) = \nFuk{1}(\rT^{*}[G/G]).
  \end{equation*}
  The heuristic we shall now appeal to is Abouzaid's theorem \cite{abouzaid-fiber-generates,abouzaid-based-loops}, which states that
  \begin{equation*}
    \Fuk(\rT^{*}Q) \cong \Mod \rC_{*}(\Omega Q),
  \end{equation*}
  (at least when $Q$ is spin). Thus we expect
  \begin{equation*}
    \nFuk{1}(\rT^{*}[G/G]) \cong \Mod \rC_{*}^{G}(\Omega G).
  \end{equation*}
  The algebra we want is actually the the Hochschild cohomology of this category, or the symplectic cohomology of $\rT^{*}[G/G]$. The heuristic for that is given by a related theorem of Viterbo and Abbondandolo-Schwarz that states
  \begin{equation*}
    \SH^{*}(\rT^{*}Q) \cong \rH_{*}(\cL Q),
  \end{equation*}
  where now $\cL$ denotes the topological loop space. Thus we expect
  \begin{equation*}
    \nHom{3}(\bone_{\Delta},\bone_{\Delta}) \cong \rH_{*}^{G}(\cL G).
  \end{equation*}
  This algebra is larger than $\rH_{*}^{G}(\Omega G)$, but since $\Omega G$ is a normal subgroup of $\cL G$, we can restrict to it to get the desired action. (This is also evident from the fact that, since $\rC_{*}^{G}(\Omega G)$ is commutative, it is contained in its derived center, that is, the Hochschild cohomology of its category of modules.)
\end{enumerate}

\begin{remark}
  The algebra $\rH_{*}^{G}(\cL G)$ is somewhat larger than $\rH_{*}^{G}(\Omega G)$. There is a semidirect product splitting $\cL G \cong \Omega G \rtimes G$, where the two factors are invariant under action of $G$ by conjugation. Therefore $\rH_{*}^{G}(\cL G) \cong \rH_{*}^{G}(\Omega G) \otimes \rH_{*}^{G}(G)$. While the action of $\rH_{*}^{G}(\Omega G)$ on $\QH^{*}_{G}(Y)$ is quantum in nature, the action of $\rH_{*}^{G}(G)$ on $\QH^{*}_{G}(Y)$ is induced by the action of $G$ on the underlying space $Y$.
\end{remark}

  

\subsection{Intrinsic mirror symmetry}

Teleman \cite{teleman-icm} proposed that the $2$-category containing gauged Fukaya categories of Hamiltonian $G$-manifolds (and more generally categories with topological $G$-action) should be equivalent to a 3D B-model category on $\BFM(G^{\vee})$, such as the Kapustin-Rozansky-Saulina $2$-category.

Now $\rH_{*}^{G}(\Omega G) = \cO(\BFM(G^{\vee}))$, and the arguments of the previous subsection show that when $Y$ is a Hamiltonian $G$-manifold, there is an action of this ring on $\QH_{G}^{*}(Y)$. Hence by taking $\Spec$ we obtain a subscheme $\cL_{Y} = \Spec \QH_{G}^{*}(Y) \subset \BFM(G^{\vee})$, which according to \cite[\S{}2.2]{teleman-coulomb} is a Lagrangian subvariety. It seems worth pointing out that Teleman's argument for this also uses the compatibility of $E_{2}$ and $E_{3}$ structures. Naturally, one expects that the 3D HMS equivalence carries the Hamiltonian $G$-manifold $Y$ over to an object of the 3D B-model $2$-category supported on $\cL_{Y}$.

To pursue the 3D HMS equivalence, we shall attempt to imitate what is called \emph{intrinsic mirror symmetry} \cite{gross-siebert} in the 2D setting, though at a higher categorical level. This attempt to formulate the 3D HMS correspondence will not be completely successful, but the difficulty that arises has a direct analog in the 2D setting (it has to do with the difference between the full $\SH^{*}(Y)$ and the degree-zero part $\SH^{0}(Y)$).

\subsubsection{2D case} A simplified story of 2D intrinsic mirror symmetry goes as follows. Suppose that $Y$ is a log Calabi-Yau manifold, with an SYZ fibration $\pi : Y \to B$, and a Lagrangian section $L \subset Y$. In this situation, we expect the Fukaya category $\Fuk(Y)$ to be monoidal, with $L$ as the monoidal unit object \cite{subotic,pascaleff,abouzaid-bottman-niu}, whose wrapped Floer cohomology is concentrated in degree zero. Thus the ``mirror algebra'' $A = \HF^{0}(L,L)$ is a commutative ring. The mirror transform arises by taking hom with $L$:
\begin{equation}
  \label{eq:intrinsic-hms}
  \HF^{*}(L,-) : \Fuk(Y) \to \Coh(\Spec(A)).
\end{equation}
If \eqref{eq:intrinsic-hms} is an equivalence we say intrinsic HMS holds; we shall assume this henceforth. For example, this is true when $Y = (\bC^{\times})^{n} \cong \rT^{*}T^{n}$, $L$ is a cotangent fiber, and $A$ is the ring of Laurent polynomials.

The ring $A$ can also be found in other places. Consider the symplectic cohomology of $Y$, which is isomorphic to the wrapped Floer cohomology of the diagonal in $\overline{Y} \times Y$, and also to the Hochschild cohomology of the Fukaya category. Define $A^{\star}$ to be this ring,
\begin{equation*}
  A^{\star} = \SH^{*}(Y) \cong \HF^{*}_{\overline{Y}\times Y}(\Delta,\Delta) \cong \HH^{*}(\Fuk(Y)).
\end{equation*}
The ring $A$ sits inside $A^{\star}$ as the degree zero part, $A = \SH^{0}(Y)$. There is also a functor
\begin{equation*}
  \Fuk(Y) \to \Mod A^{\star},
\end{equation*}
which has no reason to be an equivalence. Indeed, this is just an instance of the functor
\begin{equation*}
  \cC \to \Mod \HH^{*}(\cC)
\end{equation*}
that exists for any $k$-linear $A_{\infty}$-category $\cC$. From this perspective, the magic of mirror symmetry is that we can find a smaller algebra $A \subset \HH^{*}(\cC)$ such that restricting to this algebra makes the functor into an equivalence, and that, in this sense, the whole category can be reconstructed from a piece of its Hochschild cohomology.

\subsubsection{3D case}
\label{sec:3d-intrinsic}
Now $Y$ becomes $X = [\frakg^{*}/G]$, $\Fuk(Y)$ becomes $\nFuk{2}(X)$, and the structures are categorified. In this case we have $A^{\star}$ is the Hochschild cohomology category, and we conjecturally have
\begin{equation*}
  A^{\star} = \HH(\nFuk{2}(X)) \cong \Hom_{\overline{X} \times X}(\Delta,\Delta) = \nFuk{1}(\rT^{*}[G/G]) \cong \Mod \rC_{*}^{G}(\Omega G).
\end{equation*}
Now $\Mod \rC_{*}^{G}(\Omega G) \cong \Qcoh(\BFM(G^{\vee}))$, and by $1$-affineness of $\BFM(G^{\vee})$, modules for this category are sheaves of categories over this space, thus we obtain a functor 
\begin{equation*}
  \nFuk{2}([\frakg^{*}/G]) \to \Mod (\Mod \rC_{*}^{G}(\Omega G)) \cong \shvcat(\BFM(G^{\vee})).
\end{equation*}
This functor should not be expected to be an equivalence, since it is a functor of the form $\cC \to \Mod \HH^{*}(\cC)$, where $\cC$ is a $2$-category. The algebra $A^{\star}$ is too big. After all, the 3D B-model is not $\shvcat(\BFM(G^{\vee}))$, but it should contain things like sheaves of categories supported on Lagrangians in $\BFM(G^{\vee})$ \cite{teleman-icm}.

The case where $G = T$ is a torus is easier to describe, since $\BFM(T^{\vee}) = \rT^{*}T^{\vee}_{\bC}$. If we restrict from $\Qcoh(\rT^{*}T^{\vee}_{\bC})$ to $\Qcoh(T^{\vee}_{\bC})$ (analogous to the restriction from $\SH^{*}(Y)$ to $\SH^{0}(Y)$), we obtain a functor
\begin{equation*}
  \nFuk{2}([\mathfrak{t}^{*}/T]) \to \shvcat(T^{\vee}_{\bC}),
\end{equation*}
which is closer to what we want, but as described in \cite[Section 3.2]{teleman-icm}, this only recovers the part of the $2$-category that lives in a formal neighborhood of the zero section $T^{\vee}_{\bC} \subset \rT^{*}T^{\vee}_{\bC}$.

We have now pushed our ``intrinsic'' approach about as far as it can go without further input, so we should now consider other proposed definitions of the 3D B-model of a holomorphic symplectic manifold. As mentioned several times, Teleman \cite{teleman-icm} uses the Kapustin-Rozansky-Saulina $2$-category for this purpose. For a holomorphic symplectic target space of the form $\rT^{*}X$, one can suppose that the category has a definition directly in terms of $X$. One such proposal is to take a category of ``Ind-coherent sheaves of categories'' on $X$, notated $\indcohshvcat(X)$. This notion is to be contrasted to ``quasi-coherent sheaves of categories'' $\shvcat(X)$, a more well-understood notion \cite{gaitsgory}. A precise definition $\indcohshvcat(X)$ has been given in some cases by Arinkin, and this notion has been adopted in \cite{gammage-hilburn-mazel-gee} as the model for 3D B-branes.

In the present case $G = T$, this suggests the following conjecture.
\begin{conjecture}
  \label{conj:indcohshvcat}
  There is an equivalence of $2$-categories
  \begin{equation*}
    \nFuk{2}([\mathfrak{t}^{*}/T]) \to \indcohshvcat(T^{\vee}_{\bC}),
  \end{equation*}
  and the Hochschild cohomology $1$-categories of these $2$-categories are equivalent to $\Qcoh(\BFM(T^{\vee}))$.
\end{conjecture}


\begin{remark} 
  We shall perform a sanity check on this conjecture by heuristically calculating the Hochschild homology $1$-category of the right-hand side. Suppose that $\indcohshvcat(T^{\vee}_{\bC})$ is the $2$-category of Ind-coherent sheaves of categories on $T^{\vee}_{\bC}$. Then the Hochschild cohomology category should be the category of Ind-coherent sheaves on the derived loop space $LT^{\vee}_{\bC} \cong \rT[-1]T^{\vee}_{\bC}$,
  \begin{equation*}
    \HH(\indcohshvcat(T^{\vee}_{\bC})) \cong \Indcoh(LT^{\vee}_{\bC}) \cong \Indcoh(\rT[-1]T^{\vee}_{\bC}).
  \end{equation*}
  The fiberwise Koszul dual of $\rT[-1]T^{\vee}_{\bC}$ is $\rT^{*}[2]T^{\vee}_{\bC}$, and we have an equivalence
  \begin{equation*}
    \Indcoh(\rT[-1]T^{\vee}_{\bC}) \cong \Qcoh(\rT^{*}[2]T^{\vee}_{\bC}) \cong \Qcoh(\BFM(T^{\vee})).
  \end{equation*}
  Note that the distinction between $\Indcoh$ and $\Qcoh$ is crucial for this equivalence. This matches our claim that
  \begin{equation*}
    \HH(\nFuk{2}([\frakg^{*}/G])) \cong \Qcoh(\BFM(G^{\vee})).
  \end{equation*}
  Indeed, the arguments in this remark support this claim in a way that is independent of the Floer theory heuristics used previously.
\end{remark}


\begin{remark}
The $2$-truncated version of $\nSymp{1} = \nFuk{3}(\pt[2])$ that considers $1$-shifted symplectic manifolds as forming a $2$-category has been used by Crooks-Mayrand \cite{crooks-mayrand} to construct Moore-Tachikawa TFTs \cite{moore-tachikawa}. There is also work of Cazassus \cite{cazassus} in a similar direction. 
\end{remark}

\section{Examples of 3D HMS}
In this section we collect some more examples that show how $\nFuk{2}$ fits into the story of 3D Homological Mirror Symmetry. If we take $\nFuk{2}([\mathfrak{g}^*/G])$ as a model for the category of 3D A-branes  on  $[\mathfrak{g}^*/G]$, we can reinterpret  \cite{teleman-icm} as a study of a kind of 3D HMS for the mirror pair
\begin{equation*}
  [\mathfrak{g}^*/G] \longleftrightarrow \BFM(G^{\vee}).
\end{equation*}
Since Teleman does not explicitly use a 3D A-model, the connection between his work and 3D HMS might not be immediately apparent upon reading \cite{teleman-icm}. In this respect, as a candidate for the 3D A-model, we believe $\nFuk{2}([\mathfrak{g}^*/G])$ provides a helpful missing piece. We hope that the $2$-Fukaya category, besides being an interesting object on its own, can bring greater conceptual clarity to existing constructions and results in the literature.

We will revisit some basic elements of Teleman's picture from this perspective, focusing on the abelian case where  $G=T$ is a torus. To be clear, the examples we consider are already covered in some form by \cite{teleman-icm}, though without mention of $\nFuk{2}$. 
\cite{teleman-icm} pushes the story much further, including a beautiful discussion of Peterson-Rietsch mirror symmetry for flag varieties. Here we shall not go into these deeper aspects, as our purpose is merely to demonstrate the usefulness of thinking in terms of $\nFuk{2}$.

As we explained, the space $\BFM(T^{\vee})$ is isomorphic to the holomorphic cotangent bundle $\rT^* T^\vee_\bC$. Following \cite{teleman-icm} we model the category of 3D B-branes in $\rT^* T^\vee_\bC$ via the Kapustin-Rozansky-Saulina 2-category, notated  $\KRS(\rT^* T^\vee_\bC)$. 
As a first approximation, objects of $\KRS(\rT^* T^\vee_\bC)$ are holomorphic Lagrangian subvarieties in $\rT^* T^\vee_\bC$. They should come with extra data, such as a sheaf of categories, amounting to a 3D brane structure; for simplicity we shall omit these data in the present discussion. The Hom-category between  two Lagrangians $\mathcal{L}_1$ and $\mathcal{L}_2$ should be the category of \emph{Matrix Factorizations} (MF) associated with the derived intersection $\mathcal{L}_1 \times_{\rT^* T^\vee_\bC} \mathcal{L}_2$. 

Let us explain this point in some more detail.  The derived intersection $\mathcal{L}_1 \times_{\rT^* T^\vee_\bC} \mathcal{L}_2$ is $(-1)$-shifted symplectic, thus it can be written locally as a derived critical locus \cite{brav2019darboux}.\footnote{We should remark that, although most of this document deals with shifted symplectic structures in the \emph{differentiable} setting, here $\mathcal{L}_1 \times_{\rT^* T^\vee_\bC} \mathcal{L}_2$ is as an actual  algebraic shifted symplectic stack in the sense of \cite{ptvv}.}
Matrix factorization categories are naturally invariants of derived critical loci. Further, it had been long expected that, via a globalization step, it should be possible to associate a MF-type category to any $(-1)$-shifted symplectic stack. This was recently proved by Hennion-Holstein-Robalo \cite{hennion2024gluing}. Their result provides a  justification to the definition of the Hom-categories in the KRS 2-category. 

We can attempt to formulate 3D HMS as an equivalence of categories
\begin{equation*}
  \nFuk{2}([\mathfrak{t}^*/T]) \cong \KRS(\rT^* T^\vee_\bC),
\end{equation*}
but as it has been said, even the full and precise definition of these categories is beyond the reach of current understanding (though our Conjecture \ref{conj:indcohshvcat} is one attempt to make it more precise). However pieces of these categories are easier to describe. We shall consider the completion at the zero section in $\rT^{*}T^{\vee}_{\bC}$, and then look at some test Lagrangians in $[\mathfrak{t}^{*}/T]$ and their expected mirror partners.


\subsection{Completion at the zero section}
While 2-category $\KRS(\rT^* T^\vee_\bC)$ is a complicated object, we can obtain a simpler 2-category by completing $\rT^* T^\vee_\bC$ at the zero section. As described in \cite[Section 3.2]{teleman-icm}, this has the effect of \emph{linearizing} $\KRS(\rT^* T^\vee_\bC)$, in the sense of turning it into modules over a certain monoidal category.   
 \begin{expectation}
   \label{claimcompletion}
   The completion of $\KRS(\rT^* T^\vee_\bC)$ at the zero-section is equivalent to the 2-category 
   \begin{equation*}
     \shvcat(T^\vee_\bC) \cong \Mod (\Qcoh(T^\vee_\bC )).
   \end{equation*}
   (The equivalence of the displayed categories is the fact that $T^\vee_\bC$ is 1-affine.)
 \end{expectation}
 This completion should be compatible with the functor $\nFuk{2}([\mathfrak{t}^{*}/T]) \to \shvcat(T^{\vee}_{\bC})$ obtained in Section \ref{sec:3d-intrinsic}. 

\subsection{Test Lagrangians}
 What we shall now do is to explain how we expect this equivalence to work on some selected objects, and check that the corresponding Hom-categories match. We focus on the following test Lagrangians: 
\begin{itemize}
\item Let $F:\mathfrak{t}^* \to [\mathfrak{t}^*/T]$ be the quotient map. This is shifted Lagrangian, and coincides with the shifted cotangent fiber  $\rT^*_e[1]BT$.
\item Let $x$ be in $\mathfrak{t}^*$. Let $$L_x: BT \to [\mathfrak{t}^*/T]$$ be the map induced by the inclusion $i_x: * \to \mathfrak{t}^*$. Note that $L_0$ is the zero section of the shifted cotangent stack $[\mathfrak{t}^*/T]$. In general, $L_x$ is a translation of the zero section in the cotangent direction.  
\end{itemize}
Let $\eta: \mathfrak{t}^* \to T^\vee_\bC$ be the composition of the exponential map $\exp : \mathfrak{t}^{*} \to T^{\vee}$ and the inclusion $T^{\vee} \to T^{\vee}_{\bC}$

\begin{expectation}
\label{claimtestlag}
The shifted Lagrangian $F$ is mirror to the zero section $T^\vee_\bC \to \rT^*T^\vee_\bC$, and the shifted  Lagrangian $L_x$  is mirror to  the cotangent fiber of $\rT^*T^\vee_\bC$ at $\eta(x)$. That is,
  \begin{equation*}
    F \longleftrightarrow T^\vee_\bC, \quad  L_x \longleftrightarrow \rT^*_{\eta(x)}T^\vee_\bC.
  \end{equation*}
\end{expectation}

\begin{expectation}
As mentioned in Section \ref{sec:ham-g-man}, one can associate to any Hamiltonian $T$-manifold $M$ a shifted Lagrangian $[M/T] \to [\mathfrak{t}^{*}/T]$. Also, one can associate a holomorphic Lagrangian $\cL_M$ in $\rT^*T^\vee_\bC$, as explained in \cite{teleman-icm} and Section \ref{teleact}. This object has several different presentations:
\begin{enumerate}
\item We can view $\cL_M$ as the spectrum of the equivariant quantum cohomology of $M$:
  \begin{equation*}
    \cL_M \cong \Spec (\QH_{T}^{*}(M)).
  \end{equation*}
\item Assume that $M$ is actually a toric variety for $T$, with open orbit $U \cong T_\bC \subset M$. The mirror of $M$ is a superpotential $W:T^\vee_\bC \to \bC$. Then we can view $\cL_M$ as the graph of $dW$:
  \begin{equation*}
    \cL_M \cong \Gamma_{dW}.
  \end{equation*}
\end{enumerate}
Naturally, we expect these objects to correspond under mirror symmetry:
\begin{equation*}
  [M/T] \longleftrightarrow \cL_{M}.
\end{equation*}  
\end{expectation}


Based on Expectation \ref{claimtestlag}, we expect to have the following commutative diagrams of 2-categories
\begin{equation}
\label{diag1}
\begin{gathered}
\xymatrix{
\nFuk{2}([\mathfrak{t}^*/T]) \ar[rrrr]^-{\Hom_{\nFuk{2}([\mathfrak{t}^*/T])}(F, -)}  \ar@{-->}[d]  &&&&  \Prl \ar@{->}[d]^=  \\
\KRS(\rT^*T^\vee_\bC) \ar@{->}[rrrr]^-{\Hom_{\KRS(\rT^*T^\vee_\bC)}(T^\vee_\bC, -) } &&&& \Prl } 
\end{gathered}
\end{equation}

\begin{equation}
\label{diag2}
\begin{gathered}
\xymatrix{
\nFuk{2}([\mathfrak{t}^*/T]) \ar[rrrr]^-{\Hom_{\nFuk{2}([\mathfrak{t}^*/T])}(L_x, -)}  \ar@{-->}[d]  &&&&  \Prl \ar@{->}[d]^=  \\
\KRS(\rT^*T^\vee_\bC) \ar@{->}[rrrr]^-{\Hom_{\KRS(\rT^*T^\vee_\bC)}(\rT^*_{\eta(x)}T^\vee_\bC, -) } &&&& \Prl }
\end{gathered}
\end{equation}
where the dashed vertical arrows indicate the conjectural equivalence between the 3D A and B-models, and $\Prl$ is the $2$-category of presentable $\infty$-categories.
 We will discuss the commutativity of these diagrams   when evaluated on a toric Hamiltonian $T$-manifold $M$.

\subsubsection{The commutativity of \eqref{diag1}} We begin with the equivalence mentioned in Section \ref{sec:ham-g-man},
\begin{equation*}
  \Hom_{\nFuk{2}([\mathfrak{t}^*/T])}(F,   [M/T]) \cong \Fuk(M).
\end{equation*}
The Hamiltonian $T$-action on $M$ induces a topological $T$-action on $\Fuk(M)$. This was conjectured in \cite{teleman-icm}, and proved rigorously in \cite{savelyev-global-fuk} and \cite{oh2019continuous}. Next, as observed by Teleman, the datum of a topological  $T$-action is the same as the datum of a tensor action of $\Qcoh(T^\vee_\bC)$ on $\Fuk(M)$, where $\Qcoh(T^\vee_\bC)$ is equipped with its ordinary tensor product. By Expectation\ref{claimcompletion},  the KRS category of the completion of $\rT^*T^\vee_\bC$ at the zero section is $\Mod(\Qcoh(T^\vee_\bC))$. Thus $\Fuk(M)$ equipped with the action of $\Qcoh(T^\vee_\bC)$ defines an object inside the completed KRS category.

A key insight from  \cite{teleman-icm} is that mapping $M$ to $\Fuk(M)$ with its $\Qcoh(T^\vee_\bC)$-module structure loses information and does not capture the full 3D HMS story. Rather, as we discussed, the $T$-manifold $M$ is mirror to a Lagrangian brane $\mathcal{L}_M$ in $\KRS(\rT^*T^\vee_\bC)$ that may intersect other branes that are far away from the zero section. Only the restriction of $\mathcal{L}_M$ to  the completion $(\widehat{\rT^*T^\vee_\bC})_{T^\vee_\bC}$  is captured by the category 
\begin{equation*}
 \Fuk(M) \in \Mod(\Qcoh(T^\vee_\bC)).   
\end{equation*}
However, since we are trying to compute the Hom with the zero section inside $\KRS(\rT^*T^\vee_\bC)$, we can restrict first to a formal neighborhood of it. That is,  we can harmlessly work inside the simpler $2$-category $\Mod(\Qcoh(T^\vee_\bC))$. Once we pass to $\Mod(\Qcoh(T^\vee_\bC))$, taking the Hom with the zero section is the same as forgetting the action: it is the forgetful functor
\begin{equation*}
  \Mod(\Qcoh(T^\vee_\bC)) \to \Prl.
\end{equation*}
Plugging $[M/T] \in \nFuk{2}([\mathfrak{t}^*/T])$ into the top arrow of diagram \eqref{diag1}, we obtain $\Fuk(M)$, while if we plug $\cL_{M}$ into the lower arrow, we obtain the matrix factorization category of the superpotential $W: T^{\vee}_{\bC} \to \bC$. These two categories are equivalent by 2D HMS for toric varieties.

 \begin{example}
 Let us discuss in more detail the example of $M=S^2=\bP^1$ with its Hamiltonian $T=S^1$-action.  
 The mirror LG model of 
 $\bP^1$ is $W=x + \frac{1}{x}:\bC^* \to \bC$. The mirror holomorphic Lagrangian is
 \begin{equation*}
   \cL_{S^2}:=\Gamma_{dW}=\{y = 1 - \frac{1}{x^2}\} \subset \rT^* T^\vee_\bC   
 \end{equation*}
The Lagrangian $\cL_{S^2}$ and the zero-section $T^\vee_\bC \subset \rT^* T^\vee_\bC$ intersect transversely in the two points 
$1, -1 \in \bC^*$, and thus 
\begin{equation}
\label{homs2}
\Hom_{\KRS(\rT^* T^\vee_\bC)}(T^\vee_\bC, \cL_{S^2}) \cong \Perf(k) \oplus \Perf(k) \cong \Fuk(S^2)
\end{equation}
as expected. 

We can get to the same result arguing in a slightly different way. Since we are computing the Hom with the zero section, we can pass   to the completion $(\widehat{\rT^*T^\vee_\bC})_{T^\vee_\bC}$. This is what we did before, in the context our general discussion of the commutativity of \eqref{diag1}. This means  considering 
$\Fuk(S^2)$ equipped with its topological $S^1$-action.    Taking (twisted) fixed points for the $S^1$-action on $\Fuk(S^2)$ gives a spectral decomposition of the category. Here is how this works. Topological $S^1$-actions on $\Perf(k)$, where $k=\bC$, are classified by a point $\lambda \in \bC^*$: indeed giving a $S^1$-action is the same as giving an invertible element in $\HH^0(\Perf(k)) = \bC^*$. Let us call $\Perf(k)_\lambda$  the category $\Perf(k)$ with the $S^1$-action classified by $\lambda$. 

After completing, the intersections with the zero section can be computed by testing against cotangent fibers. That is, assuming that the Lagrangian $\cL_{S^2}$ corresponding to $\Fuk(S^2)$ is transverse to $T^\vee_\bC$ (which is the case here), we have
\begin{equation*}
  \lambda \in \cL_{S^2} \cap T^\vee_\bC \quad \text{if and only if}  \quad (\Fuk(S^2) \otimes \Perf(k)_{\lambda^{-1}})^{S^1} \neq 0.
\end{equation*}

It is not difficult to prove that, as an equivariant category, $\mathrm{Fuk}(S^2)$ decomposes as 
\begin{equation*}
  \Fuk(S^2)  \simeq \Perf(k)_{1} \oplus \Perf(k)_{-1},  
\end{equation*}
which again matches the calculation we performed with a different method in  \eqref{homs2}. 
 \end{example}
 
 \subsubsection{The commutativity of \eqref{diag2}} It turns out that, though in a different formulation, the commutativity of \eqref{diag2} is one of the main conjectures of \cite{teleman-icm}.
As mentioned in Section \ref{sec:ham-g-man}, we have  
\begin{equation*}
  \Hom_{\nFuk{2}([\mathfrak{t}^*/T])}(L_x, [M/T]) \cong \Fuk(M/\!\!/_{x}T),
\end{equation*}
while
\begin{equation*}
\Hom_{\KRS(\rT^* T^\vee_\bC)}(\rT_{\eta(x)}^{*}T^{\vee}_{\bC}, \cL_{M})   
\end{equation*}
is the ``multiplicity category'' for $\Perf(k)_{\eta(x)}$ in $\cL_{M}$. Thus the claim that diagram \eqref{diag2} commutes is equivalent to Conjecture 4.2 in \cite{teleman-icm}. That conjecture can therefore be understood as a compatibility statement in 3D HMS for $[\mathfrak{t}^*/T]$. We believe that conceptually something is gained by this reformulation, 
in that it clarifies the way this fits in the broader  picture of 3D HMS.\footnote{
There are some caveats to these last statements, already pointed out in \cite{teleman-icm}.
Note that 
\begin{equation}
\label{hom1}
\Hom_{\nFuk{2}([\mathfrak{t}^*/T])}(L_x, [M/T]) = \Fuk(M/\!\!/_{x}T)
\end{equation}
will be nonzero  only if $x$ lies in the image of the moment map $M \to \mathfrak{t}^*$. 
However  
$\Hom_{\KRS(\rT^*T^\vee_\bC)}(\rT^*_{\eta(x)}T^\vee_\bC, \cL)$
can be non-zero for general  $x$. Thus the Hom-categories
cannot be expected to match, unless $x$ is a value of the moment map. Teleman proposes to get around this issue by restricting to the \emph{unitary mirror}.}

\section{Coda: symplectic groupoids}

Let $\cG = G_{1} \rightrightarrows G_{0}$ be a Lie groupoid. When the space $G_{1}$ is equipped with a symplectic form that is multiplicative, we call $\cG$ a symplectic groupoid. According to Safronov \cite[Proposition 3.6]{safronov-lectures}, the set of symplectic groupoid structures on $\cG$ is equivalent to the space of pairs consisting of a $1$-shifted symplectic structure on the stack $[G_{0}/G_{1}]$ and a $1$-shifted Lagrangian structure on the morphism $G_{0} \to [G_{0}/G_{1}]$. From the latter data, the symplectic manifold $G_{1}$ is recovered as the derived fiber product
\begin{equation*}
  G_{1} = G_{0} \times_{[G_{0}/G_{1}]} G_{0}.
\end{equation*}
As shown by Weinstein \cite{weinstein-groupoid}, these data also induce a Poisson structure on $G_{0}$.

The main examples considered in this paper have this form.
\begin{enumerate}
\item The shifted cotangent bundle $\rT^{*}[1]M$ corresponds to the groupoid $\rT^{*}M \rightrightarrows M$, where both source and target maps are the projection. This is a bundle of abelian Lie groups, and the induced Poisson structure on $M$ is zero.
\item The coadjoint stack $[\frakg^{*}/G]$ corresponds to the action groupoid of the coadjoint action $G \times \frakg^{*} \rightrightarrows \frakg^{*}$. The induced Poisson structure on $\frakg^{*}$ is the famous Kirillov-Kostant-Souriau structure.
\end{enumerate}

In an earlier paper \cite{pascaleff}, the first author considered the formal Fukaya categorical structures that follow from a symplectic groupoid structure. The scope of the theory presented in this paper completely subsumes the results of that earlier work.

\bibliographystyle{alpha}
\bibliography{shifted}

\end{document}